\documentclass[a4paper,12pt,reqno]{amsart}

\usepackage{amsmath}
\usepackage{amssymb}
\usepackage{amsfonts}
\usepackage{graphicx}
\usepackage{mathtools}
\usepackage[colorlinks]{hyperref}
\renewcommand\eqref[1]{(\ref{#1})} 
\usepackage[pagewise]{lineno}
\graphicspath{ {images/} }
\title[Integral Rellich type inequalities]{One-dimensional integral Rellich type inequalities}
\author[T. Ozawa]{Tohru Ozawa}
\address{
	Tohru Ozawa:
	\endgraf
	Department of Applied Physics
	\endgraf
	Waseda University, Tokyo 169-8555
	\endgraf
	Japan
	\endgraf
	{\it E-mail address:} {\rm txozawa@waseda.jp}}

\author[P. Roychowdhury]{Prasun Roychowdhury}
\address{
	Prasun Roychowdhury:
	\endgraf
	Dipartimento di Matematica e Applicazioni
    \endgraf
    Universit\`a degli Studi di Milano–Bicocca
	\endgraf
	Via Cozzi 55, 20125 Milano
	\endgraf
	Italy
	\endgraf
	{\it E-mail address:} {\rm prasunroychowdhury1994@gmail.com}}

\author[D. Suragan]{Durvudkhan Suragan}
\address{
	Durvudkhan Suragan:
	\endgraf
	Department of Mathematics
	\endgraf
	Nazarbayev University
	\endgraf
	Kazakhstan
	\endgraf
	{\it E-mail address:} {\rm durvudkhan.suragan@nu.edu.kz}}

\date{\today}

\theoremstyle{plain}

\newtheorem{theorem}{Theorem}[section]

\newtheorem{lemma}{Lemma}[section]

\numberwithin{equation}{section} \allowdisplaybreaks

\usepackage[text={6in,8.6in},centering]{geometry}
\newcommand{\dr}{\:{\rm d}r}
\newcommand{\dt}{\:{\rm d}t}
\newcommand{\dta}{\:{\rm d}\tau}

\begin{document}
\begin{abstract}
The motive of this note is twofold. Inspired by the recent development of a new kind of Hardy inequality, here we discuss the corresponding Hardy--Rellich and Rellich inequality versions in the integral form. The obtained sharp Hardy--Rellich type inequality improves the previously known result. Meanwhile, the established sharp Rellich type integral inequality seems new.
\end{abstract}
\maketitle


\section{Introduction}
In the celebrated paper, \cite{hardy}, Godfrey H. Hardy first stated the famous inequality. The result reads as follows. For any $1<p<\infty$ and $f$ be a $p$-integrable function on $(0, \infty)$, which vanishes at zero, then the function $r\longmapsto \frac{1}{r}\int_{0}^{r} f(t) \dt$ is $p$-integrable over $(0, \infty)$ and there holds
\begin{align}\label{hardy}
    \int_{0}^{\infty}\bigg|\frac{1}{r}\int_{0}^{r} f(t)\dt\bigg|^p\dr\leq  \bigg(\frac{p}{p-1}\bigg)^p \int_{0}^{\infty}|f(r)|^p\dr.
\end{align}
The constant on the right-hand side of \eqref{hardy} is sharp. The development of the famous Hardy inequality \eqref{hardy} during the period 1906--1928 has its own history and we refer to \cite{kuf} (also, see the preface of \cite{hardy-book-rs}). Recent progress by Frank--Laptev--Weidl \cite{nh} presents a novel one-dimensional inequality with the same sharp constant, which improves the classical Hardy inequality  \eqref{hardy}. 

This new version looks as follows. For any $1<p<\infty$ and for any $f\in L^p(0,\infty)$, which vanishes at zero, there holds
\begin{align}\label{new-hardy}
   \int_{0}^{\infty}\sup_{0<s<\infty}\bigg|\min\biggr\{\frac{1}{r},\frac{1}{s}\biggr\}\int_{0}^{s}f(t)\dt \bigg|^p\dr\leq \bigg(\frac{p}{p-1}\bigg)^p \int_{0}^{\infty}|f(r)|^p\dr. 
\end{align}

Certainly, \eqref{new-hardy} gives an improvement of \eqref{hardy}. Recently, the multidimensional version in the supercritical case and the discrete version of \eqref{new-hardy} have been established in \cite{nhm} and \cite{nhd}, respectively. In the same spirit, one may ask about the possible structure of Hardy--Rellich and Rellich type inequalities. In this short note, we obtain the possible form of these two types of inequalities.

Let us recall the one-dimensional Hardy--Rellich inequality. For $f\in C^2(0,\infty)$ with $f^\prime(0)=0$, there holds 
\begin{align}\label{hardy-rel}
    \int_{0}^{\infty}\frac{|f^\prime(r)|^2}{r^{2}}\dr \leq 4 \int_{0}^{\infty}|f^{\prime\prime}(r)|^2\dr.
\end{align}
Starting from it, there have been several articles in which the authors studied many improvements in inequality \eqref{hardy-rel}. Here we mention only a few of them \cite{bgr-22, coss, cazacu, new-1, new-2, new-3, Owen, ya, tz} and references therein. 

Now let us write \eqref{hardy-rel} in the integral form. Note that it can be derived from the weighted one-dimensional classical Hardy inequality. This reads as follows. Let $f\in C^1(0,\infty)$, then there holds
\begin{align}\label{hardy-rel-int}
    \int_{0}^{\infty}\frac{|\int_{0}^{r}f^\prime(t)\dt|^2}{r^{2}}\dr\leq 4\int_{0}^{\infty}|f^{\prime}(r)|^2\dr.
\end{align}
Here the constant $4$  is sharp. We give an improved version of this inequality in Theorem \ref{hardy-rel-th}.

Let us briefly mention another important function inequality the so-called Rellich inequality which was first introduced in \cite{rel}. It is worth recalling the one-dimensional Rellich inequality. The classical one-dimensional Rellich inequality states that for $f\in C^2(0,\infty)$ with $f(0)=0$ and $f^\prime(0)=0$, there holds
\begin{align}\label{rel}
    \int_{0}^{\infty}\frac{|f(r)|^2}{r^{4}}\dr\leq \frac{16}{9} \int_{0}^{\infty}|f^{\prime\prime}(r)|^2\dr.
\end{align}
Over the past few decades, there has been a constant effort to improve \eqref{rel}. Here are some closely related papers \cite{hinz, ozawa, BT, MNSS, bmo, rsadv, cm}. In this short contribution, we also obtain another type of Rellich inequality (see Theorem \ref{rel-th} with $p=2$). 
To the best of our knowledge, the most recent progress in this direction was made in \cite{bmo}. However, a one-dimensional study is still missing. Thus, trying to fill this gap is another motivation for the present paper. Taking inspiration from there we obtain the following version of Rellich inequality. For any $f\in L^2(0,\infty)$ there holds
\begin{align}\label{n-rel-2}
   &\int_{0}^{\infty}\frac{1}{r^{4}}\bigg(\int_{0}^{r}\int_{0}^{\tau}|f(t)|\dt\dta\bigg)^2\dr\nonumber\\&\leq\int_{0}^{\infty}\frac{1}{r^{4}}\bigg(\int_{0}^{r}\sup_{0<s<\infty}\min\biggl\{1,\frac{\tau}{s}\biggr\}\int_{0}^{s}|f(t)|\dt\dta\bigg)^2\dr\nonumber \\&\leq \frac{16}{9} \int_{0}^{\infty}|f(r)|^2\dr.
 \end{align}
 Moreover, we will show that the constant $16/9$ is a sharp constant. Therefore, \eqref{n-rel-2} can be compared with \eqref{rel}. \color{black} Note that we have mentioned only the $L^2(0,\infty)$ case but we will discuss the result for the general $L^p(0,\infty)$ case.

\section{Preliminaries and main results}\label{prelm}
Let us begin this section with basic facts about \emph{a symmetric decreasing rearrangement}. For more details, we refer to \cite[Section 2.1]{bs}. We denote by $f^*$ the symmetric decreasing rearrangement of $f$. It is well known that $f^*$ is a nonnegative, radially symmetric, and nonincreasing function. Irrespective of several properties of $f^*$, the useful property in our context is the equimeasurability property, i.e.
\begin{align}
	\text{vol}(\{|f|>\tau\})=\text{vol}(\{f^*>\tau\}) \text{ for all }\tau\geq 0.	
\end{align}
By using the \emph{layer cake representation} and the above property, we have the following helpful identity:
\begin{align}\label{norm-presv}
	\int_{0}^\infty|f(t)|^p\dt=\int_{0}^\infty|f^*(t)|^p\dt\:\: \text{ for all }p\geq 1.
\end{align}
Also, by the simplest rearrangement inequality, for any $s>0$ there holds
\begin{align}\label{norm-comp}
    \int_0^s|f(t)|\dt\leq \int_0^sf^*(t)\dt.
\end{align}
These relations will be very much valuable in the proofs.

Now, we are ready to state the following important observation.
\begin{lemma}\label{rel-lem-1}
	Let $f$ be a locally absolutely continuous function on $(0,\infty)$. Then for a fixed $r>0$, the following identity holds:
	\begin{align}\label{rel-eq-1}
		\sup_{0<s<\infty}\min\biggl\{1,\frac{r}{s}\biggr\}\int_{0}^{s}f^*(t)\dt=\int_{0}^{r}f^*(t)\dt,
	\end{align}
	where $f^*$ is the non-increasing rearrangement of $f$.
\end{lemma}
\begin{proof}
We wish to calculate the supremum by using the monotonicity of $f^*$. For any fixed $r>0$, we consider the following two cases:

    {\bf Case 1:} Let $0<s\leq r$. Then we obtain 
	\begin{align*}
		\min\biggl\{1,\frac{r}{s}\biggr\}\int_{0}^{s}f^*(t)\dt=\int_{0}^{s}f^*(t)\dt\leq \int_{0}^{r}f^*(t)\dt.
	\end{align*}
	
	{\bf Case 2:} Let $r\leq s<\infty$. Then we have by change of variable
	\begin{align*}
		\min\biggl\{1,\frac{r}{s}\biggr\}\int_{0}^{s}f^*(t)\dt=\frac{r}{s}\int_0^{s}f^*(t)\dt\leq\frac{r}{s}\int_{0}^{s}f^*(tr/s)\dt=\int_{0}^{r}f^*(v)\:{\rm d}v.
	\end{align*}
	
	In both cases, we get 
	\begin{align*}
		\min\biggl\{1,\frac{r}{s}\biggr\}\int_{0}^{s}f^*(t)\dt \leq \int_{0}^{r}f^*(t)\dt.
	\end{align*}
	Hence the supremum is attained at $s=r$ and we arrive at 
	\begin{align*}
		\sup_{0<s<\infty}\min\biggl\{1,\frac{r}{s}\biggr\}\int_{0}^{s}f^*(t)\dt=\int_{0}^{r}f^*(t)\dt.
	\end{align*}
\end{proof}

Now we are ready to present an improvement of \eqref{hardy-rel-int}. That is, this gives a natural improvement of the Hardy--Rellich inequality in the integral form. Below we will describe the corresponding differential form which improves the original Hardy--Rellich inequality \eqref{hardy-rel} in a simple form.

\begin{theorem}\label{hardy-rel-th}
 Let $g\in C^1(0,\infty)$, then there holds
  \begin{align}\label{eqn-hardy-rel-th}
   \int_{0}^{\infty}\sup_{0<s<\infty}\bigg|\min\biggl\{\frac{1}{r},\frac{1}{s}\biggr\}\int_{0}^{s}g^\prime(t)\dt\bigg|^2\dr\leq  4 \int_{0}^{\infty}|g^\prime(r)|^2\dr.   
  \end{align}
Moreover, the constant $4$ in the above inequality is sharp in the sense that no inequality of the form
\begin{align*}
   \int_{0}^{\infty}\sup_{0<s<\infty}\bigg|\min\biggl\{\frac{1}{r},\frac{1}{s}\biggr\}\int_{0}^{s}g^\prime(t)\dt\bigg|^2\dr\leq  C \int_{0}^{\infty}|g^\prime(r)|^2\dr   
  \end{align*}
  holds, for $g\in C^1(0,\infty)$, when $C<4$.
\end{theorem}

Now we are going to discuss the second main result of this note. Before presenting the statement first let us recall the classical one-dimensional $L^{p}$-Rellich inequality. This reads as follows. Let $p>1$, $f\in C^2(0,\infty)$ with $f(0)=0$ and $f^\prime(0)=0$ there holds
\begin{align}\label{p-rel}
    \int_{0}^{\infty}\frac{|f(r)|^p}{r^{2p}}\dr\leq \frac{p^{2p}}{(p-1)^p(2p-1)^p} \int_{0}^{\infty}|f^{\prime\prime}(r)|^p\dr.
\end{align}
Now we are ready to demonstrate the one-dimensional Rellich-type inequality in the following integral form.
\begin{theorem}\label{rel-th}
 Let $f\in L^p(0,\infty)$, $p>1$. Then we have
 \begin{align}\label{eqn-rel-th}
   &\int_{0}^{\infty}\frac{1}{r^{2p}}\bigg(\int_{0}^{r}\int_{0}^{\tau}|f(t)|\dt\dta\bigg)^p\dr\nonumber\\&\leq\int_{0}^{\infty}\frac{1}{r^{2p}}\bigg(\int_{0}^{r}\sup_{0<s<\infty}\min\biggl\{1,\frac{\tau}{s}\biggr\}\int_{0}^{s}|f(t)|\dt\dta\bigg)^p\dr\nonumber \\&\leq \frac{p^{2p}}{(p-1)^p(2p-1)^p} \int_{0}^{\infty}|f(r)|^p\dr.
 \end{align}
Moreover, the constant $\frac{p^{2p}}{(p-1)^p(2p-1)^p}$ in the above inequality turns out to be sharp in the sense that no inequality of the form 
 \begin{align*}
   &\int_{0}^{\infty}\frac{1}{r^{2p}}\bigg(\int_{0}^{r}\int_{0}^{\tau}|f(t)|\dt\dta\bigg)^p\dr \leq C \int_{0}^{\infty}|f(r)|^p\dr.
 \end{align*}
 holds, for $f\in L^p(0,\infty)$ when $C<\frac{p^{2p}}{(p-1)^p(2p-1)^p}$.
\end{theorem}

\section{Proofs of Theorems \ref{hardy-rel-th} and \ref{rel-th}}\label{proof}
This section is concerned with the proofs of Theorems \ref{hardy-rel-th} and \ref{rel-th}. Before going further let us recall the following lemma.
\begin{lemma}\cite[Lemma 3.1]{nhm}\label{hr-lemma}
    Let $1<p<\infty$. Let $w$ be any nonnegative function on $(0,\infty)$. Assume $h$ is a strictly positive non-decreasing function on $(0,\infty)$ such that  $s h(r)\leq r h(s)$ for any $r,s\in(0,\infty)$ with $r\leq s$. Let $f$ be a locally absolutely continuous function on $(0,\infty)$. Then we have
		\begin{equation*}
			\int_{0}^{\infty}w(r)\sup_{0<s<\infty}\bigg|\min\biggl\{\frac{1}{h(r)},\frac{1}{h(s)}\biggr\}\int_{0}^{s}f(t)\dt\bigg|^p\dr \leq \int_{0}^{\infty}w(r)\bigg|\frac{1}{h(r)}\int_{0}^{r}f^*(t)\dt\bigg|^p\dr,
		\end{equation*}
		where $f^*$ is the non-increasing rearrangement of $f$.
\end{lemma}
Now, as a direct consequence of Lemma \eqref{hr-lemma}, we derive the proof of Theorem \ref{hardy-rel-th}.

{\bf Proof of Theorem \ref{hardy-rel-th}:} 
Let us consider $w(r)=1$ and $h(r)=r$ to be functions on $(0,\infty)$ and substitute these in Lemma \ref{hr-lemma} with $p=2$, then we have
  \begin{equation*}
		\int_{0}^{\infty}\sup_{0<s<\infty}\bigg|\min\biggl\{\frac{1}{r},\frac{1}{s}\biggr\}\int_{0}^{s}g^{\prime}(t)\dt\bigg|^2\dr \leq\int_{0}^{\infty}\frac{1}{r^{2}}\bigg|\int_{0}^{r}(g^{\prime})^*(t)\dt\bigg|^2\dr.
	\end{equation*}
	By using the Hardy--Rellich inequality in the form \eqref{hardy-rel-int} for the function $(g^\prime)^*$, we obtain
	\begin{align*}
	  \int_{0}^{\infty}\sup_{0<s<\infty}\bigg|\min\biggl\{\frac{1}{r},\frac{1}{s}\biggr\}\int_{0}^{s}g^\prime(t)\dt\bigg|^2\dr&\leq 4\int_{0}^{\infty}|{(g^\prime)}^{*}(r)|^2\dr\\&=4\int_{0}^{\infty}|g^\prime(r)|^2\dr.
	\end{align*}
In the last step, we have used norm preserving property \eqref{norm-presv}. The sharpness follows from the optimality of the constant in \eqref{hardy-rel-int}. It ends the proof.

{\bf Proof of Theorem \ref{rel-th}:}
The first inequality follows from the property of the supremum. Now taking the integral of \eqref{rel-eq-1} from $0$ to $r$  we have
\begin{align}\label{rel-eq-2}
		\int_{0}^{r}\sup_{0<s<\infty}\min\biggl\{1,\frac{\tau}{s}\biggr\}\int_{0}^{s}f^*(t)\dt\dta=\int_{0}^{r}\int_{0}^{\tau}f^*(t)\dt\dta.
\end{align}
Then
\begin{align*}
		&\int_{0}^{\infty}\frac{1}{r^{2p}}\bigg(\int_{0}^{r}\sup_{0<s<\infty}\min\biggl\{1,\frac{\tau}{s}\biggr\}\int_{0}^{s}|f(t)|\dt\dta\bigg)^p\dr\\&\overset{\eqref{norm-comp}}{\leq} \int_{0}^{\infty}\frac{1}{r^{2p}}\bigg(\int_{0}^{r}\sup_{0<s<\infty}\min\biggl\{1,\frac{\tau}{s}\biggr\}\int_{0}^{s}f^*(t)\dt\dta\bigg)^p\dr\\&\overset{\eqref{rel-eq-2}}{=}\int_{0}^{\infty}\frac{1}{r^{2p}}\bigg(\int_{0}^{r}\int_{0}^{\tau}f^*(t)\dt\dta\bigg)^p\dr\\&\overset{\eqref{p-rel}}{\leq}\frac{p^{2p}}{(p-1)^p(2p-1)^p}\int_{0}^{\infty}|f^*(r)|^p\dr\\&\overset{\eqref{norm-presv}}{=}\frac{p^{2p}}{(p-1)^p(2p-1)^p}\int_{0}^{\infty}|f(r)|^p\dr.
\end{align*}
This establishes the new type of Rellich inequality in the integral form.

{\bf Optimality:} We set 
\begin{align}\label{rel-const}
    C_{p}:=\sup_{f\in L^p(0,\infty)\setminus\{0\}}\frac{\int_{0}^{\infty}\frac{1}{r^{2p}}\big(\int_{0}^{r}\int_{0}^{\tau}|f(t)|\dt\dta\big)^p\dr}{\int_{0}^{\infty}|f(r)|^p\dr}.
\end{align}
The validity of \eqref{eqn-rel-th} immediately implies 
\begin{align*}
   C_{p}\leq \frac{p^{2p}}{(p-1)^p(2p-1)^p}.
\end{align*}
So it remains to show the reverse inequality and this will be done by giving a proper minimizing sequence. We divide the proof into some steps.

{\bf Step 1.} Let us start with a cut-off function $\chi:[0,\infty)\rightarrow \mathbb{R}$ with the following properties:
\begin{itemize}
    \item[1.] $\chi(r)\in [0,1]$ for all $r\in [0,\infty)$ and $\chi$ is smooth;
    \item[2.] $\chi$ satisfies the following
\begin{equation*}
\chi(r)=
\begin{dcases}
1 & 0\leq r\leq 1, \\
0 & 2\leq r< \infty; \\
\end{dcases}
\end{equation*}
\item[3.] $\chi$ is decreasing function, i.e. $\chi^\prime(r)\leq 0$ for all $r\in (0,\infty)$. 
\end{itemize}
Now for a small $\epsilon>0$, let us define the minimizing functions $\{f_\epsilon\}$ as follows:
\begin{align*}
f_\epsilon(r):=r^{\frac{\epsilon-1}{p}}\chi(r).
\end{align*}

{\bf Step 2.} In this step we will estimate r.h.s. of \eqref{eqn-rel-th}. The denominator of \eqref{rel-const} gives
\begin{align}\label{denom}
\int_{0}^{\infty}|f_\epsilon(r)|^p\dr&=\int_0^\infty r^{\epsilon-1}\chi^p(r)\dr\nonumber
\\&=\int_0^1 r^{\epsilon-1}\dr+\int_1^2r^{\epsilon-1}\chi^p(r)\dr\nonumber
\\&=\frac{1}{\epsilon}+O(1).
\end{align}
Therefore, for a fixed positive $\epsilon$, we have $f_\epsilon\in L^p(0,\infty)$.

{\bf Step 3.} In this part we will evaluate the numerator of \eqref{rel-const}. Using the integration by parts, we have
\begin{align}\label{numeo}
 &\int_{0}^{\infty}\frac{1}{r^{2p}}\bigg(\int_{0}^{r}\int_{0}^{\tau}|f_\epsilon(t)|\dt\dta\bigg)^p\dr\nonumber
 \\&=\int_{0}^{\infty}\frac{1}{r^{2p}}\bigg(\int_{0}^{r}\int_{0}^{\tau}t^{\frac{\epsilon-1}{p}}\chi(t)\dt\dta\bigg)^p\dr\nonumber
 \\&=\bigg(\frac{p}{\epsilon-1+p}\bigg)^p\int_{0}^{\infty}\frac{1}{r^{2p}}\bigg[\int_0^r\chi(\tau)\tau^{\frac{\epsilon-1+p}{p}}\dta-\int_{0}^{r}\int_{0}^{\tau}t^{\frac{\epsilon-1+p}{p}}\chi^\prime(t)\dt\dta\bigg]^p\dr\nonumber
 \\&\geq \bigg(\frac{p}{\epsilon-1+p}\bigg)^p\int_{0}^{\infty}\frac{1}{r^{2p}}\bigg[\int_0^r\chi(\tau)\tau^{\frac{\epsilon-1+p}{p}}\dta\bigg]^p\dr\nonumber
 \\&=\bigg(\frac{p}{\epsilon-1+p}\bigg)^p\bigg(\frac{p}{\epsilon-1+2p}\bigg)^p\int_{0}^{\infty}\frac{1}{r^{2p}}\bigg[\chi(r)r^{\frac{\epsilon-1+2p}{p}}-\int_{0}^{r}\tau^{\frac{\epsilon-1+2p}{p}}\chi^\prime(\tau)\dta\bigg]^p\dr\nonumber
 \\&\geq \bigg(\frac{p}{\epsilon-1+p}\bigg)^p\bigg(\frac{p}{\epsilon-1+2p}\bigg)^p\int_{0}^{\infty}r^{\epsilon-1}\chi^p(r)\dr\nonumber
 \\&=\bigg(\frac{p}{\epsilon-1+p}\bigg)^p\bigg(\frac{p}{\epsilon-1+2p}\bigg)^p\bigg[\int_0^1 r^{\epsilon-1}\dr+\int_1^2r^{\epsilon-1}\chi^p(r)\dr\bigg]\nonumber
 \\&=\frac{1}{\epsilon}\bigg(\frac{p}{\epsilon-1+p}\bigg)^p\bigg(\frac{p}{\epsilon-1+2p}\bigg)^p+O(1).
\end{align}
In between, exploiting $\chi^\prime\leq 0$, we used a simple inequality $(a+b)^p\geq a^p$ twice, for nonnegative real numbers $a$ and $b$ with $p>1$.

{\bf Step 4.} Finally, by using \eqref{denom} and \eqref{numeo} we deduce
the ratio
\begin{align*}
    &\frac{\int_{0}^{\infty}\frac{1}{r^{2p}}\big(\int_{0}^{r}\int_{0}^{\tau}|f(t)|\dt\dta\big)^p\dr}{\int_{0}^{\infty}|f(r)|^p\dr}\\
    &\geq \frac{\frac{1}{\epsilon}\big(\frac{p}{\epsilon-1+p}\big)^p\big(\frac{p}{\epsilon-1+2p}\big)^p+O(1)}{\frac{1}{\epsilon}+O(1)}\rightarrow \frac{p^{2p}}{(p-1)^p(2p-1)^p}\; \text{ for }\epsilon\rightarrow 0.
\end{align*}
Hence $\{f_\epsilon\}$ is the required minimizing sequence and, in turn, we have $$C_p=\frac{p^{2p}}{(p-1)^p(2p-1)^p}.$$

\medskip

\section*{Acknowledgments}
This work is partially supported by the NU program 20122022CRP1601. The first author is supported in part by JSPS Kakenhi 18KK0073, 19H00644. The second author is partially supported by National Theoretical Science Research Center Operational Plan (Project number: 112L104040).

\end{document}